\input amstex.tex 
\documentstyle{amsppt} 
\magnification=1200 
\baselineskip 16pt plus 2pt 
\advance\vsize -2\baselineskip
\parskip 2pt 
\NoBlackBoxes

\topmatter 
 
\title Removing metric anomalies from Ray Singer torsion  \endtitle 
 
\thanks   Supported in part by NSF
\endthanks

\author D. Burghelea (Ohio State University)
\endauthor 
\leftheadtext {Removing metric anomalies ...} 
\rightheadtext {D. Burghelea} 
  
\abstract 

Ray Singer torsion is a numerical invariant associated with a compact 
Riemannian manifold equipped with a flat bundle and a 
Hermitian structure on this bundle.
In this note we show how one can remove the dependence on the
Riemannian metric and on the Hermitian structure 
with the help of a base 
point and of an Euler structure, in order to obtain a topological invariant. 
A numerical invariant for an Euler structure  and  additional data is also 
constructed. 

\endabstract 
\toc 
\widestnumber\subhead{3.1} 
\head 0. Introduction \endhead 
\head 1. Euler structures \endhead 
\head 2. The Invariant $I(M,...)$ \endhead 
\head 3. Proof of the Main Theorem \endhead
\head 4. References \endhead
\endtoc

\endtopmatter 
 
\document

\proclaim{0. Introduction}\endproclaim

In this Note all vector spaces $V$ will be over the field $\Bbb R$ of real 
numbers and the scalar products are always positive definite. If $V$ has 
dimension one then the scalar products on $V$ can be identified to the 
elements of the
quotient space $V/\{-1,+1\}$ of $V$ by the antipodal action.
 
Let $(M^n,m_0)$ be a base pointed closed smooth manifold with fundamental 
group 
$\Gamma= \pi_1(M,m_0),$ and $\rho: \Gamma \to GL(V)$ a finite dimensional
representation. 

Denote by
$\Cal F_\rho = (\Cal E_\rho, \nabla_\rho)$ the associated vector bundle $\Cal E_\rho \to M$ equipped
with the canonical  flat connection whose holonomy representation at $m_0$ is given by $\rho.$ 
The fiber $\Cal E_{m_0}$ of $\Cal E$ above $m_0$ identifies canonically to 
$V.$ Note that the data $(M,m_0,\rho)$ and $(\Cal F =(\Cal E\to M, \nabla), 
m_0)$ with $\nabla$ a flat connection in the bundle $\Cal E \to M$ are 
equivalent.

Denote by $H^*(M; \rho)= H^*(M; \Cal F_\rho)$ 
the singular cohomology of $M$ with coefficients in $\rho$ and by 
$\det( H(M;\rho))$ the one dimensional vector space
$$\det( H(M;\rho)):= \bigotimes _{i=1}^n (\bigwedge^{b_i}H^i(M;\rho))^{\epsilon (i)},$$
where $b_i= \dim H^i(M;\rho)$ and $W^{\epsilon (i)}= W$ resp. $W^*,$ if $i$ 
is even resp. odd.
  
To a Riemannian metric $g$ on $M$ and a Hermitian 
structure ( smooth fiberwise scalar product) $\mu$ on $\Cal E_\rho$
one can associate:

(a) The scalar product  
$\Bbb T_{met}(M,\rho, g,\mu)$ on 
$\det H(M;\rho)$, obtained by  
Hodge theory. Indeed the Hodge theory identifies the singular cohomology 
vector spaces 
$H^q(M;\Cal F)$  to the space of harmonic q-forms. The scalar 
product 
provided by $g$ and $\mu$ on $\Omega^q(M;\rho)$ restricted to harmonic forms 
induces a scalar product on each $H^q(M;\Cal F),$ and therefore on 
$\det H(M;\rho).$ 

(b) A positive real number $T_{an}(M\rho, g, \mu),$ the Ray Singer torsion, 
obtained 
from the zeta regularized determinants of the Laplacians $\Delta_q$ 
induced by the pair 
$(g,\mu).$

The metric anomaly for $T_{an},$ cf [BZ] section iv, permits to 
show that 
in the case of odd dimensional manifolds, the scalar product 
$(T_{an}(M\rho, g, \mu))^{-1}\cdot \Bbb T_{met} (M,\rho, g,\mu)$
on the vector space $\det H(M;\rho),$ 
referred to as Ray Singer metric, is independent of $g$ and $\mu,$
hence is a topological invariant \footnote {Here topological invariant 
actually means differential invariant;  i.e. an object which although 
defined with the help of geometric data like Riemannian metric, Hermitian 
structure, etc.
is independent of them. It is 
however possible to show  that this particular invariant is in fact 
a topological invariant in the traditional sense}. 
This is not the case for an even dimensional manifold. Note that 
the Euler Poincar\'e 
characteristic of an odd dimensional manifold is zero.

In [Tu] Turaev has introduced the concept of Euler structure for compact 
smooth manifolds with 
zero Euler Poincar\'e characteristic, and has verified that $H_1(M;\Bbb Z)$ 
acts freely and  
transitively on the set $E(M)$ of Euler structures. He has shown that Euler 
structures can be also defined 
combinatorially, for 
triangulated compact spaces, and for a compact smooth manifold the analytic 
Euler structures identify canonically to the combinatorial
Euler structures. 

For a triangulated 
space $(X,\tau)$ with Euler
Poincar\'e characteristic zero, an Euler structure $A\in E(X,\tau)$ 
and a
finite dimensional representation $\rho$ over an arbitrary field,  
Farber and 
Turaev have defined in [FT] an element in $ \det H(M;\rho))/\{+1,-1\}.$
and 
have verified that it 
is a combinatorial invariant, i.e independent of triangulation up to a 
PL (piecewise linear) equivalence. 
As indicated above in the case the field is $\Bbb R,$ such  element 
is a scalar product.
They have also obtained an analytic interpretation of this element in 
the case the manifold is of odd dimension and the field is 
$\Bbb R.$ cf [FT].

Consider the one dimensional vector space $\det V^{-\chi(M)}=
(\bigwedge^{\dim V} V)^{-\chi(M)},$  ($V^{-1}= V^*,$ the dual of $V$) 
equipped with the scalar product $(\det\mu_0)^{-\chi(M)}$ induced from 
$\mu_0.$
If $\chi (M)= 0$ then $\det V^{-\chi(M)}= \Bbb R$ and the induced 
scalar product is the canonical scalar product on $\Bbb R,$
hence independent on $\mu_0.$
In this note: 

(1) We observe that Euler structures can be defined for any 
base
  pointed closed smooth manifold
$(M,m_0)$ and, as in [Tu], that $H_1(M;\Bbb Z)$ acts freely and transitively 
on the set $E(M,m_0)$ of Euler structures of $(M,m_0).$ 
If $\chi (M)=0$ the choice of 
$m_0$ is superfluous, the set $E(M,m_0)$ being canonically identified 
to the set $E(M)$ introduced by Turaev. This set consists of the equivalence classes of nonsingular vector fields modulo
the equivalence relation introduced by Turaev in [Tu].

(2) Inspired by [BZ], we define for
 an  Euler structure $A\in E(M,m_0),$ 
a  Riemannian metric $g$ and 
a Hermitian structure $\mu$ which is parallel in an unspecified  small
neighborhood of $m_0,$
a numerical invariant, $I(M,\rho, A, g,\mu)\in \Bbb R,$ 
and we verify that:

\proclaim {Main Theorem}

(1): If $M$ is odd dimensional then $I(M,\rho, A, g,\mu)$ is independent 
of $g$ and $\mu.$

(2): If the Hermitian structure $\mu$ induces a parallel Hermitian structure 
in the vector bundle 
$\det \Cal E_\rho \to M,$ then 
$I(M,\rho, A, g,\mu)= 0.$

(3): The scalar product on $\det H(M;\rho)\otimes (\det V)^{-\chi(M)}$,
defined by
$$e^{I(M,\rho, A, g, \mu)}\cdot T^{-1}_{an}(M\rho, g, \mu)\cdot(\Bbb T_{met} 
(M,\rho, g,\mu)\otimes (\det(\mu_0)^{-\chi(M)})$$
is independent of $g$ and $\mu.$

(4): If the  Euler Poincar\'e characteristic 
of $M$ vanishes then the scalar product defined in (3) is the same as the 
Farber Turaev torsion for a smooth triangulation (cf [FT]).

\endproclaim

Part (3) of the Main Theorem shows that the Ray-Singer metric properly 
modified with the help of an Euler structure $A\in E(M,m_0),$
defines a scalar product on $\det H(M;\rho)\otimes 
(\det V)^{-\chi(M)}.$ This scalar product, although defined with the help of a Riemannian metric and of a Hermitian structure is independent of them, and therefore is a topological invariant. 

Part (4) provides an analytic description of the Farber-Turaev torsion 
when the PL space is a triangulated smooth manifold and the field is $\Bbb R.$
 
Using the considerations in [BFK] one can routinely extend the 
Main Theorem from finite dimensional representations to representations 
$\rho :\Gamma \to GL_{\Cal A}(\Cal W),$ with
$\Cal W$ an $\Cal A$- Hilbertian module of finite type over a finite 
von Neumann algebra $\Cal A$
and $(M,\rho)$ of ``determinant class'' in the sense of [BFK]. 

I wish to thank Misha Farber for bringing to my attention the concept of Euler structure and for very helpful comments concerning this Note. 
 
\subhead {1. Euler structures }
\endsubhead

Let $(M^n,m_0)$ be a base pointed closed smooth manifold. Consider 
$\Cal V(M,m_0)$ be the collection of smooth vector fields $X$ on $M$ with $X(x)\ne 0$ for 
$x\ne m_0.$ $\Cal V(M,m_0)$ is nonempty. We say that 
$X_1, X_2 \in \Cal V(M,m_0)$
are equivalent iff there exists a smooth family of vector fields $X_t \in   \Cal V(M,m_0),$
$1\leq t\leq 2;$ equivalently iff there exists an embedded disc $D^n\subset M^n$ 
centered at $m_0$
and a continuous  family of smooth vector fields $X_t \in   \Cal V(M,m_0),$
$1\leq t\leq 2,$ so that $X_t(x) \ne 0$ for $x\in M\setminus \text{Int} D^n.$ 
Denote the collection of equivalence classes of such vector 
fields by $E(M,m_0).$ 


There exists a natural action
$$\lambda: H_1(M;\Bbb Z)\times E(M,m_0) \to E(M,m_0)$$
noticed by Turaev in the case $\chi(M)=0$ which is free and transitive.
To describe this action we proceed as follows:
Choose a (generalized) triangulation $\tau$ of $M$ 
in the sense of [BFKM] section 5, $\tau= (h,g')$ with $h$ a selfindexing Morse 
function and $g'$ a Riemannian metric so that the vector field 
$-grad_{g'}h$ satisfies the 
Morse Smale condition, (cf. [BFKM]).The cells =simplexes of this 
triangulation are the unstable manifolds of the critical points of 
$h$ with respect to the vector field $ -grad_{g'}h.$ It is known that
any simplicial triangulation can be obtained in this way.

An ``embedded spray'' will be given by a smooth 
embedding  $\alpha: K \to M$, of the base point union $K= 
\bigvee _{y\in Cr(h)} I_y ,$ 
($I_y$ denotes the 
based pointed interval $([0,1], 0)$) with $\alpha_y:= \alpha: I_y \to M$ 
a smooth embedding so that  $\alpha_y(0)=m_0$ and 
$\alpha_y(1)= y.$
Clearly a smooth regular neighborhood of $\alpha(K)$ is 
an embedded disc $D^n\subset M^n.$

To a given embedded spray $\alpha$ one can associate  
an Euler structure $A(\alpha)\in E(M,m_0)$ defined as follows: First one 
chooses an 
embedded disc $D^n$ (smooth image of the unit disc in $\Bbb R^n$),
which 
is a regular neighborhood of the  smooth spray $\alpha (K).$ Second  one
defines the vector field $X$ on $M\setminus D^n$ by $X= -grad_g(h).$ Third 
one extends $X$ radially and scaled by
a factor equal to the distance 
to the center inside $D^n.$  The resulting vector field $X$ will have only 
one zero at the center of the disc, hence it represents an Euler structure. 
Different regular neighborhoods 
lead to  vector fields representing the same Euler structure,
and so do  embedded sprays which are isotopic by isotopy 
fixed on the vertices of $K$ and $m_0.$
By changing the spray one changes 
the Euler structure. 
In fact for any given Euler structure one can choose an embedded spray so 
that the 
vector field constructed above represents the given Euler structure.
If $\alpha$ and $\beta$ are two embedded sprays, the singular cycle 
$\sum_{y\in Cr(h)}(-1)^{index (y)}(\alpha_y  -\beta_y)$ represents a homology class 
$u(\alpha, \beta) \in H_1(M;\Bbb Z)$ and given 
$u \in H_1(M;\Bbb Z)$ and $\alpha$ an embedded spray, 
there exists an embedded spray $\beta$ so that $u= u(\alpha, \beta).$ 
We put
$$\lambda (u(\alpha, \beta), A(\alpha))= A(\beta).$$
The verification that any Euler structure is of the form $A(\alpha)$ and that 
$\lambda$ is a well defined action, free and transitive, can be 
done as in [Tu]. It can be also derived by using 
elementary obstruction theory\footnote {choose a generalized triangulation $\tau$ and assume 
that $m_0$
lies in a top dimensional
simplex (cell). Denote by $M(k)$ the $k-$th skeleton.There is no obstruction 
to construct the homotopy 
$X_t$ between $X_1$ and $X_2,$  ($X_1, X_2 \in \Cal V(M,m_0)$)
above an open neighborhood of 
$M(n-2)$ with $X_t(x) \ne 0$ for any $t\in [1,2].$ The obstruction to 
extend such homotopy 
to $M(n-1)$ is an $(n-1)-$cocycle in 
$C^*(M,\tau; \Cal O)$ the geometric complex provided by $\tau.$ The cocycle
is cohomological to zero iff 
the restriction of $X_t$ to $M(n-3)$ can be extended to a 
neighborhood of $M(n-1).$ Fixing a vector field $X_1\in \Cal V(M,m_0)$
one can provide $X_2$ and homotopy $X_t$ realizing such cycle 
(cohomology class) as obstruction. 
Once such homotopy is constructed one can easily show that $X_1$ and 
$X_2$ define the same Euler class.}.

As already stated embedded sprays which are isotopic by isotopy 
which are fixed on the vertices of $K$ and $m_0$ give rise to the same Euler 
structure. Consequently one can consider (when $dim M>2$) continuous sprays  
by replacing 
smooth embedding of $I$ by continuous maps ( since a continuous 
map of $I$ can be 
approximated arbitrary close in $C^0$ topology by a smooth embedding unique 
up to an isotopy). In particular if $\gamma: [0,1] \to M$
with $\gamma(0)= \gamma(1)$ represents the element $[\gamma]\in H_1(M;\Bbb Z)$
and $\alpha^{\gamma}$ denotes the continuous spray defined by 
$\alpha^{\gamma}_y= \gamma(2t)$ for $0\leq t\leq 1/2$ and $\alpha_y(2t-1)$ for 
$1/2\leq t \leq 1,$ one can verify (cf [FT]) that 
$$\leqalignno{A(\alpha^{\gamma})= \lambda(\chi(M)[\gamma], A(\alpha)). 
&&(1.1)\cr}$$
It is not hard to see that starting with an embedded spray $\alpha$ one can produce the embeded spray $\alpha^{\gamma}$ which satisfy (1.1)even in the case 
$dim M=2.$ 

\subhead {2.  The Invariant  $ I(M,\rho, A, g,\mu)$}
\endsubhead 

Let $\pi: E\to M$ be a vector bundle
of rank $k$ and denote by 
$\Cal O$ the orientation bundle, which is in a natural way a flat bundle 
of rank one and by 
$\pi^*(\Cal O)$
its pullback on $E.$  For forms $\omega$ with values in $\Cal O$ resp. 
$\pi^*(\Cal O),$ denote by $d\omega$ the covariant differentiation with 
respect to the natural flat connection of $\Cal O$ resp. $\pi^*(\Cal O).$
Consider pairs $\tilde{\nabla}= (\nabla, \mu)$ consisting of a connection 
$\nabla$ and a $\nabla-$parallel hermitian structure $\mu$ in 
$\pi: E\to M.$ Recall that a Hermitian structure means a scalar product in each fiber $E_x$, $x\in M$ depending smoothly on $x.$
For any such $\tilde{\nabla}$ 
consider:

(1) the Euler integrand $e(\tilde \nabla) \in \Omega^k(M;\Cal O),$ which is 
a closed form,

(2) Mathai-Quillen form 
$ \Psi (\tilde \nabla) \in \Omega^{k-1}(E\setminus 0;\pi^*(\Cal O)),$ cf [BZ] 
page 
40-44.)
\newline If $\tilde \nabla_1$ and $\tilde \nabla_2 $ are two such pairs 
consider 

(3) $c(\tilde \nabla_1,\tilde \nabla_2) \in  
\Omega^{k-1}(M; \Cal O)/d (\Omega^{k-2}(M;\Cal O))$ the Chern Simon class. 

\noindent It is well known that 
$$\leqalignno{d c(\tilde \nabla_1,\tilde \nabla_2)
= e(\tilde \nabla_1) -e(\tilde \nabla_2) &&(2.1)\cr}$$ and that
if there exists 
$\phi(t):E\to E,$ $1\leq t\leq 2$ a smooth family of bundle isomorphisms 
so that $\phi(1)= Id$ and $\phi(2)^*(\tilde \nabla_1)= \tilde \nabla_2$
then $$\leqalignno{c(\tilde \nabla_1,\tilde \nabla_2)=0. &&(2.2)\cr}$$ Here 
$\phi(2)^*(\tilde \nabla_1)$ denotes the action of $\phi(2)$ on the pair 
$\tilde \nabla_1.$

If $X: M \to E\setminus 0$ is a smooth map 
denote by $X^*(\Psi (\tilde \nabla))
\in \Omega^{k-1}(M;\Cal O)$ the pullback of 
$\Psi (\tilde \nabla)$ by $X.$
We have (cf [BZ] pages 40-44) 
$$\leqalignno{ d \Psi (\tilde \nabla) = \pi^* e(\tilde \nabla) &&(2.3)\cr}$$
$$\leqalignno{ \Psi (\tilde \nabla_1)-  \Psi (\tilde \nabla_2) = 
\pi^*  c(\tilde \nabla_1,\tilde \nabla_2), &&(2.4)\cr}$$
with $\pi^* c(\tilde \nabla_1,\tilde \nabla_2)  \in 
\Omega^{k-1}(E\setminus 0;\pi^*(\Cal O))/d (\Omega^{k-2}
(E\setminus 0;\pi^*(\Cal O)).$ 
In particular we have 
$$\leqalignno{ d X^* (\Psi (\tilde \nabla)) = e(\tilde \nabla)&&(2.3')\cr}$$
and
$$\leqalignno{ X^*(\Psi (\tilde \nabla_1))-  X^*(\Psi (\tilde \nabla_2)) =  
c(\tilde \nabla_1, \tilde \nabla_2)
\in \Omega^{k-1}(M;\Cal O)/d (\Omega^{k-2}(M;\Cal O))&&(2.4')\cr}$$

If $f:M\to \Bbb R_+$ is a smooth function then
$$\leqalignno{ (fX)^*(\Psi(\tilde \nabla))= X^*(\Psi(\tilde \nabla)). &&(2.5)
\cr}$$

Suppose $X_1$ and $X_2$ are two smooth maps as above and suppose
that there exists  
$\phi(t) : E\to E,$  $1\leq t \leq 2,$ a smooth family of bundle isomorphisms
so that $\phi(1)= Id$ and $\phi (2)\cdot X_1 =X_2.$  
Then
$$\leqalignno{  X_1^*(\Psi (\tilde \nabla))-  X_2^*(\Psi (\tilde \nabla)) =
X_1^*(\Psi (\tilde \nabla)) - X_1^*(\phi(2)^*\Psi(\tilde \nabla)) 
= &&(2.6)\cr}$$
$$= c(\tilde \nabla, 
\phi(2)^*\tilde \nabla) = 0.$$
The vanishing of $c(\tilde \nabla, \phi(2)^*\tilde \nabla)$ 
follows from the definition of secondary 
characteristic classes.

These considerations will be applied to $\pi:E\to M$ the tangent bundle of 
\newline $M\setminus D^n(1/2) $ where $D^n$ will be an embedded disc centered 
at $m_0,$
$D^n(1/2)$ the disc of radius half of the radius of $D^n $
and $\tilde \nabla (g):= (\nabla (g), g),$ where $\nabla (g)$ is the 
Levi- Civita connection associated with a Riemannian metric $g$ and the Hermitian structure is given by $g.$

Let $(\Cal E \to M, \nabla)$ be a flat bundle. Given two Hermitian structures $\mu_1$ and $\mu_2$ 
one defines the smooth real valued function 
$\log\det(\mu_1, \mu_2) \in \Omega^0(M)$ by the formula 
$$\leqalignno{\log\det(\mu_1, \mu_2)(x):= 
\log\text{Vol}(Id: (\Cal E_x,\mu_1(x))\to (\Cal E_x, \mu_2(x))) &&(2.7)\cr}$$
where $ \text{Vol}(Id: (\Cal E_x,\mu_1(x))\to (\Cal E_x, \mu_2(x)))$ is the 
$\mu_2(x)-$volume of a parallelepiped generated by a $\mu_1(x)-$ 
orthonormal base. 
Equivalently $ \log \text{Vol}(Id: (\Cal E_x,\mu_1(x))\to 
(\Cal E_x, \mu_2(x)))=
1/2 \log \det (Id^*\cdot Id)$ where $Id^*$ is the adjoint of 
$Id: (\Cal E_x,\mu_1(x))\to (\Cal E_x, \mu_2(x)).$ 

Given a Hermitian structure $\mu$ one defines the closed form $\theta(\mu) 
\in \Omega^1(M),$ first introduced by Kamber-Tondeur (cf [BZ]),
as follows:
For any $x\in M$ choose $U$ a contractible open neighborhood and 
denote by $\tilde{\mu}_x$ the Hermitian structure
in  $(\Cal E|_U  \to U, \nabla|_U)$ obtained by parallel transport of the scalar product $\mu_x$ in $\Cal E_x$
with respect to the flat connection $\nabla|_U.$ Define $\theta(\mu) |_U$ by
$$\leqalignno{\theta(\mu) |_U:= d \log\det(\mu, \tilde{\mu}_x). &&(2.8)\cr}$$
One can verify that the definition is independent of the choice of $U,$ 
and that it leads to a globally defined closed one
form. It follows from (2.8) that 
$\theta(\mu) $ vanishes on open sets $V\subset M$
if $\mu|_V$ is parallel and 
$$\leqalignno{\theta(\mu_1) - \theta(\mu_2)= d\log \det (\mu_1,\mu_2). 
&&(2.9)\cr}$$

Consider now a Riemannian metric $g$ on $M$ and a Hermitian structure $\mu$ in $\Cal E_\rho$ 
which is parallel
in an open neighborhood of $m_0.$ The form $\theta(\mu)$ vanishes in that 
neighborhood. Let 
$X\in \Cal V(M,m_0).$
The quantity 
$$\leqalignno{I(M,\rho, X, g, \mu):= 1/2 \int_{M\setminus m_0} \theta(\mu)
\wedge X^*\Psi(\tilde \nabla(g)) \in  \Bbb R &&(2.10)\cr}$$ 
is well defined since the integrand has compact support.

\proclaim {Proposition 2.1}
Let $X_1$ and $X_2$ be two vector fields in 
$\Cal V(M,m_0)$ which represent the same Euler structure $A\in E(M,m_0)$.
Then
$I(M,\rho, X_1, g, \mu)= I(M,\rho, X_2, g, \mu).$
\endproclaim

Proof: Since $I(M,\rho, X, g, \mu)$ is a continuous function in $X$ 
with respect to the $C^1-$ topology, it suffices to prove the result 
for $X_1$ and $X_2$
close enough in this topology.
In view of (2.3) one can suppose that the tangent $X_1(x)$ and $X_2(x)$ 
have the 
same length (say equal to one) for $x\in M\setminus D^n(1/2)$. 
Since we can suppose 
$X_1$ and $X_2$ are sufficiently closed one can find the bundle isomorphism 
$\phi(t) :TM|_{M\setminus d^n(1/2)} \to TM|_{M\setminus D^n(1/2)},$  $1\leq t \leq 2$ so that 
$\phi(1)=Id,$ $\phi(2)\cdot X_2= X_1,$ and 
the result follows from (2.5) and (2.6). 

In view of the above  proposition  we will write 
$I(M,\rho, A, g, \mu)$ for $I(M,\rho, X, g, \mu)$ when  $X\in A,$  
$A\in E(M,m_0).$

\proclaim {Proposition 2.2}
Suppose $g,g'$ are Riemannian metrics on $M$ and $\mu, \mu'$ are 
Hermitian structures 
parallel in some neighborhood of $m_0.$ Then
$$\leqalignno{I(M,\rho, A, g, \mu) - I(M,\rho, A, g', \mu)= 
1/2 \int_{M} \theta(\mu)\wedge c(\tilde \nabla(g),
\tilde \nabla(g'))&&(2.11)\cr}$$ 
and
$$\leqalignno{I(M,\rho, A, g, \mu)- I(M,\rho, A, g, \mu')=-1/2 
\int_M  \log\det(\mu,\mu'))\cdot e(\tilde \nabla(g)) +&&(2.12)\cr}$$
$$+ 1/2 \chi(M)\log\det (\mu,\mu')(m_0).$$
\endproclaim

Proof: (2.11) follows from (2.4'). 
To prove (2.12) we proceed as follows:
Choose $D^n$ an embedded disc centered at $m_0$ 
so that $\mu$ and $\mu'$ are parallel above $D^n,$ and choose 
a vector 
field $X\in \Cal V(M,m_0)$ representing $A.$ Denote by $D^n(1/2)$ the 
disc of radius half of 
the radius of $D^n.$
In view of the definition (2.10), and because $\mu, \mu''$ are parallel on 
$D^n,$ and because of (2.9),    
$$2(I(M,\rho, A, g, \mu)- I(M,\rho, A, g, \mu'))= 
\exp( \int_{M\setminus D^n(1/2)} (\theta(\mu) - \theta (\mu'))\wedge  
X^*(\Psi(\nabla(g))),$$
$$=\exp( \int_{M\setminus D^n(1/2)}d\log \det(\mu,\mu')
\wedge X^*(\Psi(\tilde \nabla(g))$$
which in view of Stokes Theorem and of (2.3) is equal to
$$\exp( \int_{\partial D^n(1/2)} \log\det(\mu,\mu') 
X^*(\Psi(\tilde \nabla(g))) - 
\int_{M\setminus D^n(1/2)} \log\det(\mu,\mu')  
e(\tilde \nabla(g)).$$
Because $\log\det(\mu,\mu')$ is constant on $D^n,$
the last quantity equals to 
$$\exp( \log\det(\mu,\mu')(m_0)\int_{\partial D^n(1/2)}  
X^*(\Psi(\tilde \nabla(g))) - \int_{M\setminus D^n(1/2)} \log\det(\mu,\mu')  
e(\tilde \nabla(g))),$$
which by (2.3) and by Stokes Theorem is equal to
$$\exp( \log\det(\mu,\mu')(m_0)\int_{M\setminus  D^n(1/2)}  
e(\tilde \nabla(g)) - \int_{M\setminus D^n(1/2)} \log\det(\mu,\mu')  
e(\tilde \nabla(g))).$$
The result follows by adding and subtracting inside exp(....)
the quantity $$\int_{D^n(1/2)} \log\det(\mu,\mu')  
e(\tilde \nabla(g))$$  and by using the fact that 
$\int_M e(\tilde \nabla(g))= \chi(M).$

\proclaim{3. Proof of the Main Theorem} \endproclaim

Proof of (1): Since for an odd dimensional manifold $e(\tilde \nabla(g))$ 
and $c(\tilde \nabla(g),\tilde \nabla(g'))$ are both zero, the invariant 
$I(M,\rho, A, g, \mu)$ is independent of $g$ and $\mu.$

{\bf Remark 3.1:} If $[\theta(\rho)]\in H^1(M;\Bbb R)$ is the cohomology 
class represented by the form \footnote { the notation is justified by 
(2.9) which implies that the cohomology class of $\theta(\mu)$ depends 
only on $\rho;$ this cohomology class is given by the composition of
$\log \text{det}: GL(V)\to \Bbb R$ with the 
holonomy representation $\rho:\pi_1(M,m_0)\to
GL(V).$} 
  $\theta(\mu)$ and $\gamma\in H^{n-1}(M\setminus{m_0};\Bbb R)= H^{n-1}(M; \Bbb R)$ is the Euler class 
of the bundle normal to $X\in \Cal V(M,m_0),$ for $X$ representing $A,$ then 
$I(M,\rho, A, g, \mu)$ is $[\theta(\rho)]\cup \gamma [M]$ evaluated on the fundamental class cf [FT].

Proof of (2): follows from the definition (2.10). 
 
Proof of (3): For a pair $(g, \mu)$ with $\mu$ parallel in the neighborhood of 
$m_0$ consider the scalar product 
$$T_{an}(M,\rho,g,\mu)\cdot\Bbb T_{met}(M,\rho, g,\mu).$$ Note that:

(a) The quotients for the scalar products corresponding to $(g,\mu)$
and $(g',\mu)$ is by [BZ] (page 10) and  
Proposition 2.2 (1) equal to 

$$\exp (I(M,\rho,A,g,\mu)- I(M,\rho,A,g',\mu)).$$ 
This proves the independence of $g.$ 

(b) The quotients for the scalar products corresponding to $(g,\mu)$
and $(g,\mu')$ is by [BZ] (page 10),  
$$-1/2 \exp (\int_M \log \det(\mu, \mu') e(\nabla(g)).$$
By Proposition 2.2 (2) 
the difference $I(M,\rho, A, g, \mu)- I(M,\rho, A, g, \mu')$ is equal to 
$$-1/2 
\int_M  \log\det(\mu,\mu'))\cdot e(\tilde \nabla(g)) +
1/2\chi(M)\log\det (\mu,\mu')(m_0).$$
This implies the independence of $\mu.$

Proof(4): We choose a generalized triangulation $\tau =(h,g),$ a spray
$\alpha$
representing the Euler structure $A$ with respect to $\tau,$ 
$(e_1, e_2, \cdots, e_k)$ a base of $\Cal E_{m_0}$
and $D^n$ an embedded disc centered at $m_0$ which is a regular 
neighborhood of the spray $\alpha.$
We will take as representative of $A$ a vector field $X= -grad_g(h)$
on $M\setminus D^n.$ We take as metric $g$ the one provided by $\tau$ 
and as a Hermitian structure $\mu$ one which is parallel above $D^n$ 
and at $m_0$ makes the base $(e_1, e_2, \cdots, e_k)$ orthonormal.
We will calculate our invariant (scalar product defined in the 
Main Theorem (3)) with the help of these data and compare the result 
with the Farber-Turaev torsion which we recall below.

{\bf The FT(Farber Turaev) torsion.}  
Let $(M^n,m_0)$ be a base pointed 
closed smooth manifold and $\tau =(g,h)$ a 
generalized triangulation.

Choose a base $(e_1, e_2, \cdots, e_k)$ for $\Cal E_{m_0}$
and a spray $\alpha$ representing  
the 
Euler structure
$A \in E(M,m_0)$ with respect to $\tau.$
Use the parallel 
transport to provide a 
base in 
each fiber $\Cal E_y$, $y\in Cr(h).$
Once this is done we have a base in each component of the geometric complex 
$C^*(M,\tau;\rho).$    
These bases define an element $\omega_1 \in
\det H(M;\rho) $ as described in [FT] and the base 
$(e_1, e_2, \cdots, e_k)$ gives an element $\omega_2\in \det V,$
$\omega_2=e_1\wedge e_2, \cdots\wedge e_k.$ 
The element $\omega_1 \otimes \omega_2^{-\chi(M)}\in 
(\det H(M;\rho)\otimes (\det V)^{\chi(M)} /\{-1,+1\}$
is the Farber Turaev torsion, and is independent of all choices 
but the PL structure of $M$ and the Euler structure $A\in E(M,m_0).$ 
Strictly speaking in [FT] only the 
case  $\chi(M)=0$ was considered, but the above definition is implicit in 
that paper.

We observe that $\mu_y$ is exactly the scalar product which makes orthonormal
the base obtained from 
$(e_1, e_2, \cdots, e_k)$ by 
parallel transport along $\alpha_y$ and therefore the scalar products 
induced by the bases constructed in each component of $C^*(M,\tau;\rho)$  
are the same as those determined by $\tau$ and $\mu.$
 
We use the combinatorial Laplacians associated with $\tau$ and $\mu$ 
to  identify the harmonic 
elements in $C^*(M,\tau;\rho)$ with the singular cohomology and 
to obtain  a scalar product in 
$\det H(M;\rho).$  We multiply this scalar product by  the positive real 
number $(T_{comb}(M,\tau, \rho, \mu))^{-1}$
(as defined in [BFKM]) and tensor by the scalar product induced from
$\omega_2$; we obtain  a scalar product 
on $\det H(M;\rho)\otimes V^{-\chi(M)}.$ Elementary linear algebra 
permits to show that this scalar product is the same as Farber Turaev torsion.
 
It is immediate that the quotient of the scalar product defined in (3) 
and the scalar product described above is the same as the quotient between 
Ray-Singer 
torsion and Reidemeister torsion and therefore, by [BZ], is exactly 
$ 1/2 \int_{M\setminus m_0} \theta(\mu)
\wedge X^*\Psi(\tilde \nabla(g))= I(M,\rho, A, g,\mu).$
This finishes the proof of statement (4). 

Remark 3.1 show that our analytic interpretation is consistent with the 
one given by Farber Turaev for odd 
dimensional manifolds.

\proclaim{ 6. References} \endproclaim

\enddocument